\documentclass[12pt]{amsart}
\usepackage{fullpage}
\usepackage{url}

\begin{document}

\title{A proof of Dejean's conjecture}

\author{James Currie and Narad Rampersad}

\address{Department of Mathematics and Statistics \\
University of Winnipeg \\
515 Portage Avenue \\
Winnipeg, Manitoba R3B 2E9 (Canada)}

\email{\{j.currie,n.rampersad\}@uwinnipeg.ca}

\thanks{The first author is supported by an NSERC Discovery Grant.}
\thanks{The second author is supported by an NSERC Postdoctoral
Fellowship.}


\date{\today}

\begin{abstract}
We prove Dejean's conjecture.  Specifically, we show that Dejean's
conjecture holds for the last remaining open values of $n$, namely
$15 \leq n \leq 26$.
\end{abstract}

\maketitle

\section{Introduction}

Repetitions in words have been studied since the beginning of the
previous century \cite{thueI,thue}. Recently, there has been much
interest in repetitions with fractional exponent
\cite{brandenburg,carpi,dejean,longfrac,krieger,mignosi}. For
rational $1<r\le 2$, a {\bf fractional $r$-power} is a non-empty
word $w=pe$ such that $e$ is the prefix of $p$ of length
$(r-1)|p|$. We call $e$ the {\bf excess} of the repetition.
We also say that $r$ is the {\bf exponent} of the repetition $pe$.
For example, $010$ is a $3/2$-power, with excess 0. A basic problem
is that of identifying the repetitive threshold for each alphabet
size $n>1$:

\begin{quote}
What is the infimum of $r$ such that an
infinite sequence on $n$ letters exists, not containing any factor
of exponent greater than $r$?
\end{quote}

This infimum is called the {\bf repetitive threshold} of an $n$-letter
alphabet and is denoted by $RT(n)$. Dejean's conjecture \cite{dejean}
is that
$$
RT(n)=
\begin{cases}
7/4,&n=3\\
7/5,&n=4\\
n/(n-1),&n\ne 3,4.
\end{cases}
$$

Thue, Dejean and Pansiot, respectively \cite{thue,dejean,pansiot},
established the values $RT(2)$, $RT(3)$, $RT(4)$. Moulin Ollagnier
\cite{ollagnier} verified Dejean's conjecture for $5\le n\le 11$,
and Mohammad-Noori and Currie \cite{morteza} proved the conjecture
for $12\le n\le 14$.  Recently, Carpi \cite{carpi} showed that Dejean's
conjecture holds for $n\ge 33$. The present authors strengthened Carpi's
construction to show that Dejean's conjecture holds for $n\ge
27$ \cite{currie,archiv}. In this note we show that in fact
Dejean's conjecture holds for $n\ge 2$. We will freely assume the usual notions of combinatorics on
words as set forth in, for example, \cite{acw}.

\section{Morphisms}\label{morphisms}
Given previous work, it remains only to show that  Dejean's conjecture holds for $15\le n\le 26$. This follows from the fact that the following morphisms are `convenient' in the sense of \cite{ollagnier}.
To make our exposition self-contained, we demonstrate in the remainder
of this paper how these morphisms are used to prove Dejean's conjecture
for $15\le n\le 26.$  We introduce several simplifications
and one correction to the work of Moulin Ollagnier \cite{ollagnier}.

\tiny

\begin{eqnarray*}h_{15}(0) &=& 01101101011011011011010101101010110110110110110101101101\\
h_{15}(1) &=& 10101011011011010110101101101011010110110110101101010101\\
\\
h_{16}(0) &=& 101010110110110101101101101010101010101010110101010101010101\\
h_{16}(1) &=& 011010110110110110101101101010110110101010110101010101010101\\
\\
h_{17}(0) &=& 1010101010101010101101101011011010101010101011010110110110101101\\
h_{17}(1) &=& 1010101010101010101101101101101010101010110110110110110110110110\\
\\
h_{18}(0) &=& 10101010101101101101010110110101011011011010101011010110110101010101\\
h_{18}(1) &=& 01101010101101101010110110110101011011011010101011010101010101010101\\
\\
h_{19}(0) &=& 101010101010101010101101101010110110101010101010101101011010110110101101\\
h_{19}(1) &=& 101010101010101010101101101011011010101010101011011011011010110110110110\\
\\
h_{20}(0) &=& 1010101010101010101011011011010101010101011011011011011010110101011011010101\\
h_{20}(1) &=& 1010101010101010101011011011010101101101011011011011010110110101011011010110\\
\\
h_{21}(0) &=& 101010101010101010101011011010101010101101101010101010110110110110101010101010101101\\
h_{21}(1) &=& 101010101010101010101011011010101010110110101010101010110101101010101010101010110110\\
\\
h_{22}(0) &=& 101010101010101010101011010101010101010101010101101101101101101011011011011011010101\\
h_{22}(1) &=& 101010101010101010101011010101010101011011010101101101101101011011011011011011010110\\
\\
h_{23}(0) &=& 1010101010101010101010101010101010101011011010110110110110101011010110110110110110101101\\
h_{23}(1) &=& 1010101010101010101010101010101010101101101010110110110110110110110110110110110110110110\\
\\
h_{24}(0) &=& 10101010101010101010101011010101101101010101010101011010101010110110101101101101011011010101\\
h_{24}(1) &=& 10101010101010101010101011010101101101010110110101011010101010110101101101101101011011010110\\
\\
h_{25}(0) &=& 101010101010101010101010101011011010101011011010110110110101101101011010101010101011011010110110\\
h_{25}(1) &=& 101010101010101010101010101011011010101010110110110110110101101101101101101010101011011010101101\\
\\
h_{26}(0) &=& 1010101010101010101010101011010101010101101101010101010110110110101101011010110110110110110110110101\\
h_{26}(1) &=& 1010101010101010101010101011010101010101101101101101010110110110101101010110110110110110110110110110
\end{eqnarray*}
\normalsize

We remark that the last letter of $h_n(0)$ is different from the last letter of $h_n(1)$ in each case.
We also note that for each $n$, $|h_n(0)| = 4n-4$, except for $n=21$ where we have $|h_n(0)|=4n$.

\section{Maximal repetitions}
For each $h_n$ of Section~\ref{morphisms}, word 011 is a factor of $h_n(0)$ and $110$ is a factor of $h_n(1)$. It follows that $|h_n^m(1)|$ becomes arbitrarily large as $m$ increases, and that every factor of $h_n^\omega(0)$ is a factor of $h_n^m(1)$ for some $m$.

Let an occurrence of $v$ in $h_n^\omega(0)$ be written $h_n^\omega(0) = xv{\bf y}$. Suppose that $v$ has period $q$. We can
write $x = x'x''$, ${\bf y}=y'{\bf y''}$ such that $x''vy'$ has
period $q$, and $|x^{\prime\prime}vy'|$ is maximal. This is possible since every factor
of $h_n^\omega(0)$ is a factor of $h_n^m(1)$ for some $m$, and word 1 has two distinct left
extensions 01 and 11, and two distinct right extensions 10 and 11.
We refer to $x''vy'$ as the {\bf maximal period $q$ extension} of the occurrence $xv{\bf y}$ of $v$.

\section{Pansiot encoding}
Fix $n\ge 2.$ Let $\Sigma_n = \{1, 2, \ldots, n\}$. Let
$v\in\Sigma_n^*$ have length $m\ge n-1$, and write $v=v_1v_2\cdots
v_m$, $v_i\in \Sigma_n$. In the case where every factor of $v$ of
length $n-1$ contains $n-1$ distinct letters, we define the {\bf
Pansiot encoding of $v$} to be the word $b(v) = b_1b_2\cdots
b_{m-(n-1)}$ where for $1\le i\le m-n+1$
$$b_i=\left\{\begin{array}{ll}0,&v_i=v_{i+n-1}\\1,&\mbox{otherwise. }
\end{array}\right.$$

We can recover $v$ from $b(v)$ and $v_1v_2\ldots v_{n-1}$. We see
that $v$ has period $q$ if and only if $b(v)$ does. The exponent
$|v|/q$ of $v$ corresponds to an exponent ${\displaystyle
{|v|-n+1\over q}}$ of $b(v)$.

Let $S_n$ denote the symmetric group on $\Sigma_n$ with identity
{\tt id} and left multiplication, i.e.,
$$(fg)(i) = f(g(i)) \mbox{ for }f,g\in S_n,i\in\Sigma_n.$$
Let $\sigma:\{0,1\}^*\rightarrow S_n$ be the semigroup
homomorphism generated by
\begin{eqnarray*}
\sigma(0)&=&\left(\begin{array}{cccccc} 1&2&\cdots&(n-2)&(n-1)&n\\
2&3&\cdots&(n-1)&1&n
\end{array}\right)\\
\sigma(1)&=&\left(\begin{array}{cccccc} 1&2&\cdots&(n-2)&(n-1)&n\\
2&3&\cdots&(n-1)&n&1
\end{array}\right)
\end{eqnarray*}

One proves by induction that
\begin{eqnarray}\label{kernel}
\sigma(b(v))&=&\left(\begin{array}{cccccc} 1&2&\cdots&(n-2)&(n-1)&n\\
v_{m-n+2}&v_{m-n+3}&\cdots&v_{m-1}&v_m&\hat{v}
\end{array}\right)
\end{eqnarray}
where $\hat{v}$ is the unique element of
$\Sigma \setminus \{v_m,v_{m-1},\ldots,v_{m-n+2}\}$.

Suppose that $PE\in\Sigma_n^*$ is a repetition of period $q=|P|>0$ with
$|E|\ge n-1$. It follows from (\ref{kernel}) that
$\sigma(b(P))={\tt id}$; i.e.\ that $P$ is in the kernel of
$\sigma$. We refer to $b(PE)$ as a {\bf kernel repetition} of
period $q$. Conversely, if $u\in\Sigma_n^*$ and $b(u)$ is a kernel
repetition of period $q$, then we may write $u=PE=EP'$ for some
words $P,P',E$ where $|P|=|P'|=q$.

Suppose that for a morphism $h:\{0,1\}^*\rightarrow \{0,1\}^*$ there is a $\tau\in S_n$ such that
$$
\begin{array}{lcl}\tau\cdot\sigma(h(0))\cdot\tau^{-1}&=&\sigma(0)\\
\tau\cdot\sigma(h(1))\cdot\tau^{-1}&=&\sigma(1)
\end{array}
$$
In this case we say that $h$ satisfies the `algebraic condition'.

\section{Kernel repetitions with markable excess}\label{markable}
Let a uniform morphism $h:\{0,1\}^*\rightarrow \{0,1\}^*$ be given. Let $|h(0)|=r>0$. A word $v\in\{0,1\}^*$ is {\bf markable} (with respect to $h$) if whenever $h(X)xv$ and $h(Y)yv$ are prefixes of $h^\omega(0)$ with $|x|,|y|<r$, then $x = y$. If a word is markable, its extensions are markable.
Let $U$ be the set of length 2 factors of $h^\omega(0)$.  A word $v\in\{0,1\}^*$ is {\bf 2-markable} (with respect to $h$) if whenever
\begin{enumerate}
\item $u$, $u'\in U$,
\item $h(X)xv$ is a prefix of $h(u)$ with $|x|<r$, and
\item $h(Y)yv$ is a prefix of $h(u')$ with $|y|<r$,
\end{enumerate}
then $x = y$.

If $|v|= r$ and $v$ is a factor of $h^\omega(0)$, then $v$ is a factor of $h(u)$, some $u\in U$. It follows that if $v$ is 2-markable, then $v$ is markable. For each $n$, if $h=h_n$, we find $U=\{01,10,11\}$. It follows that all length $r$ factors $v$ are factors of $h(0110)$. A finite check shows that
if $|v|= r$ and $v$ is a factor of $h^\omega(0)$, then $v$ is 2-markable, hence markable.

Let $n$ be fixed, $15\le n\le 26$ and let $h=h_n$. One checks that
$h$ satisfies the algebraic condition. Suppose that $v=pe$ is a
kernel repetition with period  $q=|p|$, where $h^\omega(0)=xv{\bf y}$.
Notice that every length $q$ factor of $pe$ is conjugate to $p$,
by the periodicity of $pe$. It follows that every length $q$
factor of $pe$ lies in the kernel of $\sigma$. Suppose that the
excess $e$ of $v$ is markable. Let $V= x^{\prime\prime}vy'$ be the
maximal period $q$ extension of the occurrence $xv{\bf y}$ of $v$. Write
$x=Xx'$, ${\bf y}=y'{\bf Y}$, so that $h^\omega(0)=XV{\bf Y}$. Write $V=PE=EP'$
where $|P|=q.$ Since $E$ is an extension of $e$, $E$ is markable.
Write $X=h(\chi)\chi'$ where $|\chi'|<r$ and write $XP =
h(\gamma)\gamma'$ where $|\gamma'|<r$. It follows from the
markability of $E$ that $\chi'=\gamma'$. Then the maximality of
$V$ yields $|\chi'|=|\gamma'|=0$. We may thus write $X=h(\chi)$,
$E=h(\eta)\eta'$, with $|\eta'|<r$. By the maximality of $V$, word
$\eta'$ must be the longest common prefix of $h(0)$ and $h(1)$.
Since $E$ is a prefix and suffix of $PE$ and $E$ is markable, we
know that $r$ divides $|P|$. In total then, we may write
$XPE=h(\chi\pi\eta)\eta'$ where $h(\pi)=P$, and $\eta$ is a prefix
of $\pi$. Also, since $h$ satisfies the algebraic condition,
$\sigma(\pi)={\tt id}$. Thus $\pi\eta$ is a kernel repetition in
$h^\omega(0)$. We see that $|PE|=r|\pi\eta|+|\eta'|$.

The maximality of $V$ implies that $\pi\eta$ is maximal with
respect to having period $|\pi|$. This means that if $\eta$ is
markable, we can repeat the foregoing construction. Eventually we
obtain a kernel repetition $\mathcal{P}\mathcal{E}$ with
non-markable excess $\mathcal{E}$. If it takes $s$ steps to arrive
at $\mathcal{P}\mathcal{E}$ then we find that
$|PE|=r^s|\mathcal{P}\mathcal{E}|+|\eta'|\sum_{i=0}^{s-1}r^{i}$
and $|P|=r^s|\mathcal{P}|$.

\section{Main result}

Let $n$ be fixed, $15 \leq n \leq 26$ and let $h = h_n$.
Suppose that $u_1$ is a factor of $h^\omega(0)$ with $|u_1|=\ell.$
Extending $u_1$ by a suffix of length at most $r-1$, and a prefix
of length at most $r-1$, we obtain a word $h(u_2)$, some factor
$u_2$ of $h^\omega(0)$, where $|u_2|\le \lfloor
(\ell+2(r-1))/r\rfloor$. Repeating the argument, we find that
$u_1$ is a factor of $h^2(u_3)$, some factor $u_3$ of
$h^\omega(0)$ where

\begin{equation}\label{iterate}
|u_3|\le \left\lfloor
\frac{\left\lfloor(\ell+2(r-1))/r\right\rfloor+2(r-1)}{r}
\right\rfloor.
\end{equation}

Define $$I(\ell,r) = \left\lfloor
\frac{\left\lfloor(\ell+2(r-1))/r\right\rfloor+2(r-1)}{r}
\right\rfloor.$$

Let ${\bf w}$ be the $\omega$-word over $\Sigma_n$ with prefix
$123\cdots (n-1)$ and Pansiot encoding $b({\bf w}) = h^\omega(0)$.
We will show that ${\bf w}$ contains no $\left({n\over
n-1}\right)^+$-powers. Suppose to the contrary that $pe$ is a
repetition in ${\bf w}$ with $|pe|/|p|
> n/(n-1)$ and $e$ a prefix of $p$.

First suppose that $|e|\geq(n-1)$.  Let $PE=b(pe)$. Then $PE$ is a
kernel repetition. Let $\eta'$ be the longest common prefix of
$h(0)$ and $h(1)$.  As in the previous section, replacing $pe$ and
$PE$ by longer repetitions of period $|P|$ if necessary, we may
assume that $h^\omega(0)$ contains a kernel repetition
$\mathcal{P}\mathcal{E}$ with non-markable excess $\mathcal{E}$
such that
$|PE|=r^s|\mathcal{P}\mathcal{E}|+|\eta'|\sum_{i=0}^{s-1}r^{i}$
and $|P|=r^s|\mathcal{P}|$.

We find that
\begin{eqnarray*}
1+{1\over n-1}&=&{n\over
n-1}\\
&<&{|pe|\over|p|}\\
&=&{|PE|+n-1\over|P|}\\
&=&{r^s|\mathcal{P}\mathcal{E}|+|\eta'|\sum_{i=0}^{s-1}r^{i}+n-1\over
r^s|\mathcal{P}|}\\
&=&{r^s|\mathcal{P}|+r^s|\mathcal{E}|\over
r^s|\mathcal{P}|}+{|\eta'|\sum_{i=1}^{s}r^{-i}\over
|\mathcal{P}|}+{n-1\over
r^s|\mathcal{P}|}\\
&<&1+{1 \over |\mathcal{P}|} \left(|\mathcal{E}|+|\eta'|{r\over
r-1}+n-1\right)
\end{eqnarray*}
so that
\begin{eqnarray*}
|\mathcal{P}|&<&(n-1)\left({|\mathcal{E}|}
+{|\eta'|} {r\over r-1}+{n-1}\right)\\
\end{eqnarray*}
and
\begin{eqnarray*}
|\mathcal{P}\mathcal{E}|&<&|\mathcal{E}|+(n-1)\left({|\mathcal{E}|}
+{|\eta'|} {r\over r-1}+{n-1}\right)\\
&\le&r+(n-1)\left(r
+(r-1) {r\over r-1}+{n-1}\right)\\
&\le&4n +(n-1)(9n-1)\\
&=&9n^2-6n+1.
\end{eqnarray*}

We use that $|\mathcal{E}|<r$ (since all factors of $h^\omega(0)$
of length $r$ or greater are markable) and $r\le 4n$ (as observed
in Section~\ref{morphisms}). Finally, since $\eta'$ is a proper
prefix of $h(0)$, $|\eta'|<r$.

 One verifies that $I(9n^2-6n+1,r)=2$. Since every length
2 factor of $h^\omega(0)$ is a factor of $0110$, word $b(PE)$ must
be a factor of $h^2(0110)$. Let $v$ be the word of $\Sigma_n$ with
prefix $123\cdots (n-1)$ and Pansiot encoding $h^2(0110)$. Since
$b(PE)$ is a kernel repetition, word $v$ contains a repetition
$\hat{p}\hat{e}$ with $|\hat{e}|\ge n-1$. However, a computer
search shows that $v$ contains no such repetition.

We conclude that $|e|\le n-2$. In this case,
\begin{eqnarray*}
{n\over n-1}<{|pe|\over |p|}
&\implies&|e|n>|pe|\\
&\implies&(n-2)n-(n-1)>|b(pe)|\\
&\implies&n^2-3n+1>|b(pe)|
\end{eqnarray*}
 However, $n^2-3n+1<9n^2-6n+1$, so that again $b(pe)$ must be a
factor of $h^2(0110)$, and $v$, defined as in the previous case,
must contain a $\left({n\over n-1}\right)^+$-power. However, a
computer search shows that word $v$ is $\left({n\over
n-1}\right)^+$-power free.

We have proved the following:

\noindent{\bf Main Result:} Let ${\bf w}$ be the word over
$\Sigma_n$ with prefix $123\cdots (n-1)$ and Pansiot encoding
$b({\bf w}) = h^\omega(0)$.  Word ${\bf w}$ contains no
$\left({n\over n-1}\right)^+$-powers.

\section{Final Remarks}
Our result builds on that of \cite{ollagnier}, but uses somewhat
simpler arguments, taking advantage of properties of our specific
morphisms. In addition, we have specified bounds for the various
computer checks, rather than invoking mere decidability.

A large simplification results from the fact that our morphisms
give binary words with no kernel repetitions at all (even of small
exponent). When moving from $PE$ to $\pi\eta$ in
Section~\ref{markable} one can give the relationship between the
exponents of these two kernel repetitions.
$${|PE|\over|P|}={|\pi\eta|\over |\pi|}+{|\eta'|\over r|\pi|}.$$
If it takes $s$ steps to arrive from repetition $PE$ to a
repetition $\pi\eta$ with non-markable excess, then the exponents
differ by $${|\eta'|\over|\pi|}\sum_{i=1}^sr^{-i}.$$ In the
notation of \cite{ollagnier}, $PE$ corresponds to
$\mu^s(\pi,\eta)$,  and has the largest exponent among the
$\mu^i(\pi,\eta)$, $0\le i\le s$. Unfortunately, \cite{ollagnier}
is marred by getting this backward, saying that for uniform
morphisms the largest exponent occurs either for $i = 0$ or for
$i=1$!

In fact, for the morphisms given for $n=5,6,7$, $\eta'$ is empty, so the aforementioned
reversal has no effect. However, for $8\le n\le 11$, $\eta'$ is
non-empty, and a more complicated check than indicated in
\cite{ollagnier} is necessary to ensure that the given
constructions work. Happily, they do indeed work, as a more
careful check shows.

Finally, we mention a few points regarding the search strategy for
finding morphisms. The second step of the strategy
indicated in \cite{ollagnier} calls for enumerating all candidate morphisms
of short enough length. A priori, this involves enumerating all  binary
words of length at most $r$ which are Pansiot encodings of $\left({n\over n-1}\right)^+$-free words over $\Sigma_n$. Initially this was
part of our strategy. Unfortunately, our experience supports the conjecture in
\cite{shur}, that the number of these words grows approximately as $1.24^r$
(independently of $n$.) 

For successive $r$ values we looked at  all possible pairs $\langle h(0), h(1)\rangle$ such that $|h(0)|,|h(1)|\le r$ where $h(0),h(1)$ were Pansiot encodings of $\left({n\over n-1}\right)^+$-free words and satisfied
 the algebraic condition; this allowed us to
verify the claim of \cite{ollagnier} that the morphisms presented therein for $5\le n\le 11$ are shortest possible `convenient morphisms'; the uniforms are all uniform, with 
lengths around $4n-4$ in each case. However,  storing all legal Pansiot encodings up to length $4n-4$
fills up a laptop with
2G RAM at around $n=15$. Therefore, our search program had to
migrate to computers with more and more RAM, simply to store Pansiot encodings. On the plus side,
we found a great number of `convenient morphisms'  for $12\le n\le 17$, not just the ones presented in this paper.

To find morphisms for $n$ up to $26$ (and indeed for various other higher values of $n$) we adopted a different strategy. Using backtracking, we found  legal
Pansiot encodings of length exactly $r=4n-4$ (or $r=4n$, in the case $n=21$), but only saved encodings $v$ for which the permutation $\sigma(v)$ was an $r$-cycle
(and thus a candidate for $h(1)$) or an $(r-1)$-cycle (and thus a candidate for $h(0))$. As soon as a candidate for $h(i)$ was found, it was tested together with 
each previously found candidate for  $h(1-i)$ to see whether a `convenient morphism' could be formed, in which case the search terminated. This search used very little memory, and 
terminated quickly. For $n=26$, our $C^{++}$ code found the morphism in just over 6 hours.
\section{Acknowledgments}

We would like to thank Dr.\ Randy Kobes for facilitating access to
computational resources.
Some of the calculations were performed on the WestGrid high performance
computing system (\url{www.westgrid.ca}).

We have recently been informed that Dr.\ Micha\"el Rao has also announced a proof
of Dejean's conjecture.

\end{document}